\documentclass[11pt,twoside]{article}

\usepackage{amsbsy,amsfonts,amsmath,amssymb,eucal,mathrsfs}
\usepackage[all]{xy}
\usepackage{pstricks}
\usepackage{pst-plot}

\newtheorem{definition}{Definition}[section]

\newtheorem{prop}[definition]{Proposition}

\newtheorem{lemm}[definition]{Lemma}

\newtheorem{coro}[definition]{Corollary}
\newtheorem{theo}[definition]{Theorem}
\newtheorem{notation}[definition]{Notation}

\newtheorem{construction}[definition]{Construction}

\newtheorem{remark}[definition]{Remark}

\newtheorem{remarks}[definition]{Remarks}

\newtheorem{example}[definition]{Example}

\newtheorem{examples}[definition]{Examples}

\newtheorem{nothing}[definition]{$\!\!$}

\newtheorem{definition*}{Definition}[section]
\newenvironment{defi*}{\begin{definition*} \rm}{\end{definition*}}
\newtheorem{definitions*}[definition*]{Definitions}
\newenvironment{defis*}{\begin{definitions*} \rm}{\end{definitions*}}
\newtheorem{prop*}[definition*]{Proposition}
\newtheorem{lemm*}[definition*]{Lemma}
\newtheorem{coro*}[definition*]{Corollary}
\newtheorem{theo*}[definition*]{Theorem}
\newtheorem{remark*}[definition*]{Remark}
\newenvironment{rema*}{\begin{remark*} \rm}{\end{remark*}}
\newtheorem{remarks*}[definition*]{Remarks}
\newenvironment{remas*}{\begin{remarks*} \rm}{\end{remarks*}}
\newtheorem{example*}[definition*]{Example}
\newenvironment{exam*}{\begin{example*} \rm}{\end{example*}}
\newtheorem{examples*}[definition*]{Examples}
\newenvironment{exams*}{\begin{examples*} \rm}{\end{examples*}}
\newtheorem{nothing*}[definition*]{$\!\!$}
\newenvironment{noth*}{\begin{nothing*} \rm}{\end{nothing*}}

\begin{document}

\def\pt{\{{\rm pt}\}}
\def\ra{\rightarrow}
\def\s{\sigma}\def\OO{\mathbb O}\def\PP{\mathbb P}\def\QQ{\mathbb Q}
 \def\CC{\mathbb C} \def\ZZ{\mathbb Z}\def\JO{{\mathcal J}_3(\OO)}
\newcommand{\G}{\mathbb{G}}
\def\proof{\noindent {\it Proof.}\;}
\def\qed{\hfill $\square$} 
\def \uh {{\widehat{u}}}
\def \vh {{\widehat{v}}}
\def \wh {{\widehat{w}}}
\def \xh {{\widehat{x}}}



\newcommand{\N}{\mathbb{N}}
\newcommand{\Z}{\mathbb Z}
\newcommand{\R}{\mathbb{R}}
\newcommand{\Q}{\mathbb{Q}}
\newcommand{\C}{\mathbb{C}}
\renewcommand{\H}{\mathbb{H}}
\renewcommand{\O}{\mathbb{O}}
\newcommand{\F}{\mathbb{F}}

\renewcommand{\a}{{\cal A}}
\newcommand{\az}{\a_\Z}
\newcommand{\ak}{\a_k}

\newcommand{\rc}{\R_\C}
\newcommand{\cc}{\C_\C}
\newcommand{\hc}{\H_\C}
\newcommand{\oc}{\O_\C}

\newcommand{\rk}{\R_k}
\newcommand{\ck}{\C_k}
\newcommand{\hk}{\H_k}
\newcommand{\ok}{\O_k}

\newcommand{\rz}{\R_Z}
\newcommand{\cz}{\C_Z}
\newcommand{\hz}{\H_Z}
\newcommand{\oz}{\O_Z}

\newcommand{\RR}{\R_R}
\newcommand{\CR}{\C_R}
\newcommand{\HR}{\H_R}
\newcommand{\OR}{\O_R}

\newcommand{\re}{\mathtt{Re}}


\newcommand{\dual}{{\bf v}}
\newcommand{\com}{\mathtt{Com}}
\newcommand{\rg}{\mathtt{rg}}

\newcommand{\g}{\mathfrak g}
\newcommand{\h}{\mathfrak h}
\renewcommand{\u}{\mathfrak u}
\newcommand{\n}{\mathfrak n}
\newcommand{\e}{\mathfrak e}
\newcommand{\plie}{\mathfrak p}
\newcommand{\q}{\mathfrak q}
\newcommand{\liesl}{\mathfrak {sl}}
\newcommand{\so}{\mathfrak {so}}

 \title{Quantum cohomology of minuscule homogeneous spaces
   II\linebreak Hidden symmetries}
 \author{P.E. Chaput, L. Manivel, N. Perrin}

\maketitle

\begin{abstract}
We prove that the quantum cohomology ring of any minuscule or
cominuscule homogeneous space, once localized at the quantum
parameter, has a non trivial involution mapping Schubert classes to
multiples of Schubert classes. This can be stated as 
a strange duality property for the Gromov-Witten invariants, 
which turn out to be very symmetric. 
\end{abstract}

 {\def\thefootnote{\relax}
 \footnote{ \hspace{-6.8mm}
 }
 }

\section{Introduction}

This paper is a sequel to \cite{cmp}, where we began a unified 
study of the quantum cohomology of (co)minuscule homogeneous
manifolds. One intriguing feature of quantum cohomology, as 
already observed by several authors and especially for 
Grassmannians (see \cite{}), 
is that it is {\it more symmetric} that ordinary cohomology. 

In this paper we show that the main tools introduced in \cite{cmp}
allow to state, and prove, a general {\it strange duality} 
statement for the quantum cohomology of a (co)minuscule homogeneous
manifold $X=G/P$. Recall that a $\ZZ$-basis for the ordinary
cohomology ring $H^*(X)$ (or for the Chow ring $A^*(X)$) of $X$ 
is given by the Schubert classes
$\sigma(w)$, where $w\in W_X$ belongs to the set of minimal lengths
representatives of $W/W_P$, the quotient of the Weyl group $W$ of $G$
by the Weyl group $W_P$ of $P$. To be precise, we let $\sigma(w)$ be 
the class of the Schubert variety $X(w_Xw)$, where $w_X$ denotes the 
longest element in $W_X$. So $\sigma(w)$ has degree $\ell(w)$ in the 
Chow ring (and twice this degree in cohomology). Note that the map
$w\mapsto p(w)=w_0w$, for $w\in W/W_P\simeq W_X$,  defines Poincar{\'e} 
duality on $X$. The map $w\mapsto \iota(w)=w_0^Pw$, where $w_0^P$
denotes the longest element of $W_P$, will play a crucial role in 
the sequel. For $w=s_{\alpha_1}\cdots s_{\alpha_{\ell(w)}}$ a reduced
  decomposition of $w\in W_X$, we let 
$$y(w)=\prod_{i=1}^{\ell(w)}n_{\alpha_i}(\alpha_0)^{\epsilon(\alpha_i)}.$$
Here $\alpha_0$ denotes the highest root of $G$, and
$n_{\alpha_i}(\alpha_0)$ is the coefficient of $\alpha_i$ when
$\alpha_0$ is written in the basis of simple roots. 
Moreover, we have let $\epsilon(\alpha)=1$
if $\alpha$ is a long root, $\epsilon(\alpha)=-1$ if $\alpha$ is short 
(in the simply-laced case, all roots are considered long). 
The rational number $y(w)$ is well 
defined since in the (co)minuscule case, reduced decompositions are
uniquely defined up to commutation relations. Finally, we let 
$\delta(w)$ be the number of occurences, in a reduced decomposition
of $w$, of the simple root $\beta$ that defines $P$. 

The Schubert classes are still a basis over $\ZZ[q]$ of the (small) 
quantum Chow ring $QA^*(X)$, whose associative product is
defined in terms of 3-points Gromov-Witten invariants. Denote by 
$QA^*(X)_{loc}$ its localization at $q$, that is,
$$QA^*(X)_{loc}=QA^*(X)\otimes_{\ZZ[q]}\ZZ[q,q^{-1}].$$
Strange duality  can be stated in a 
uniform way as follows. 

\begin{theo}[Strange Duality]\label{SD}
Let $X$ be a minuscule or cominuscule homogeneous space. 
The endomorphism $\iota$ of $QA^*(X)_{loc}$, defined  by 
$$\iota(q)=y(s_{\alpha_0})q^{-1} \quad and  \quad 
\iota(\sigma(w))=q^{-\delta(w)}y(w)\sigma(\iota(w)),$$
is a ring involution. 
\end{theo}

This result can be stated as a symmetry property of the Gromov-Witten 
invariants, see Corollary \ref{symGW}, which unexpectedly 
relates certain numbers of small degree rational curves 
with numbers of high degree rational curves. In fact we will
obtain these results by observing that the quantum product with the
class of a point maps any Schubert class to another Schubert class,
multiplied by some power of $q$, see Theorem \ref{point}. 

For Grassmannians these results were first proved by Postnikov \cite{Po};
in this case $y(w)=1$ for any $w$. For classical Grassmannians 
our strange duality statement will be deduced from the quantum 
Pieri formulas of \cite{BKT}. For the two exceptional minuscule
spaces, we have used the presentations of the quantum Chow rings
obtained in \cite{cmp}. The case of the Cayley plane has been 
checked by hand, but that of the Freudenthal variety required the 
help of a computer. 

As a consequence of this theorem, we deduce from the
formula for the smallest power of $q$ appearing in the quantum product
of two Schubert classes, 
(see \cite{FW} or \cite{cmp}), a formula for the highest power of $q$
in such a product (Corollary \ref{qmax}). For $u$ and $v$ in $W_X$, let us
denote, following W. Fulton and C. Woodward \cite{FW}, by
$\delta(u,v)$ the minimal degree of a rational curve meeting two
general translates of $X(u^*)$ and $X(v^*)$ (see also \cite{cmp} page 20
and corollary 4.12 for a combinatorial description). 

\begin{theo}\label{Qmax}
For $u,v\in W_X$, the maximal power of $q$ that appears in the 
quantum product of Schubert classes $\sigma(u)*\sigma(v)$, is
$$d_{max}(u,v)=\delta(u)-\delta(\iota(u),p(v))
=\delta(v)-\delta(\iota(v),p(u)).$$
\end{theo}

\section{A partition of the Hasse diagram}

Let $X=G/P$ be a (co)minuscule homogeneous variety. In \cite{cmp} we
defined the perimeter $d_{max}$ of $X$. For any non negative integer 
$d\leq d_{max}$, we introduced  certain 
Schubert subvarieties $T_d$ and $Y_d$ of $X$, with $T_d\subset
Y_d^*$. These varieties allowed us to define a {\it quantum Poincar{\'e}
duality} as follows: for any Schubert subvarieties $X(u), X(v)\subset
X$, the Gromov-Witten invariant $I_d(Y_d^*,X(u), X(v))$ is non 
zero if and only if $X(u), X(v)$ are contained in the (smooth) variety
$T_d$, and define Poincar{\'e} dual classes in $T_d$ -- in which case the
invariant equals one. In particular 
$$I_d(Y_d^*,\pt, T_d)=1.$$
One of the main themes of this paper is to investigate another 
quantum Poincar{\'e} type duality, defined by the non vanishing of the 
Gromov-Witten invariants $I_d(\pt,X(u),X(v))$. Let $\delta(u)$ 
denote the maximal integer $d$ such that $X(u)\subset Y_d^*$. 
The following facts are taken from \cite{cmp} page 20:
\begin{enumerate}
\item There exists a degree $d$ rational curve in $X$ joining the base
  point $e(1)=P/P$ of $X=G/P$ with the base point $e(u)=uP/P\in X(u)$, 
if and only if $d\ge\delta(u)$. 
\item If $\beta$ is the simple root defining $P$, $\delta(u)$ is the
  number of occurences of $s_{\beta}$ in a reduced decomposition of
  $u$.  
\item $\delta(u)=\delta(w_X,u)$, the minimal length of a Bruhat chain 
from $u$ to $w_X$, as defined in \cite{FW}.
\end{enumerate}

We add to this list the following essential property. 

\begin{lemm}\label{distance}
We have $T_d\subset X(w)\subset Y_d^*$ if and only if
$d=\delta(w_X,w)$. 
\end{lemm}
 
\proof 
We proved in \cite[corollary 4.12]{cmp} a combinatorial
characterisation of the smallest power appearing in a quantum
product. This gives $\delta(w_X,w)=\min\{d\ /\ T_d\subset X(w)\}$. The
lemma will follow from the fact that for all $d\leq d_{max}$, we have
$T_{d-1}\not\subset Y_d^*$ and the equivalence
$$T_{d-1}\not\subset X(w)\Leftrightarrow X(w)\subset Y_d^*.$$
These two results come from the following facts on the quivers
$Q_{Y_d^*}$ and $Q_{T_d}$ (cf. \cite{cmp} for definitions and results
on quivers):
\begin{itemize}
\item the quiver $Q_{Y_d^*}$ is obtained from the quiver $Q_X$ by
  removing all the vertices above the vertex $(\theta(\beta),d)$ where
  $\beta$ is the simple root defining $X$ and $\theta$ is the Weyl
  involution 
\item The vertices of $Q_{T_{d-1}}$ are those under the vertex
$(\theta(\beta),d)$.
\end{itemize}
The vertex $(\theta(\beta),d)\in Q_X$ is in the quiver $Q_{T_{d-1}}$ but not
in the quiver $Q_{Y_d^*}$ proving that $T_{d-1}\not\subset
Y_d^*$. Furthermore, the condition $T_{d-1}\not\subset X(w)$ is
equivalent to the fact that the vertex $(\theta(\beta),d)$ is not in the
quiver of $X(w)$ which is also equivalent to the inclusion
$X(w)\subset Y_d^*$.\qed

\medskip
A nice consequence is that we get a partition of $W_X$ in $d_{max}+1$ Bruhat
intervals,  $$W_X=\bigsqcup_{d=0}^{d_{max}}[T_d,Y_d^*].$$
We will denote by $W_d\subset W_X$ the interval $[T_d,Y_d^*]$. 
For example, $W_1$ is the image 
in $W/W_P\simeq W_X$ of the set of reflections $s_{\alpha}$, 
$\alpha$ a root of $G$. In particular $T_1$ is represented 
by $s_{\alpha_0}$. 

Recall that in \cite{cmp}, we proved 
that Gromov-Witten invariants of degreed $d$ on $X$ can be interpreted
as classical intersection numbers on an auxiliary variety $F_d$. 
This variety $F_d$ is homogeneous under the same group $G$ as $X$, 
and there is an incidence diagram 
$$\xymatrix{I_d\ar[r]^{p_d}\ar[d]^{q_d}&X\\
 F_d&}$$
Then $Z_d\subset F_d$ is the (image by $q_d$ of the) fiber of $p_d$,  
and $Y_d\subset X$ is the (image by $p_d$ of the) fiber of $q_d$. 
These varieties are given by the following table. 

$$\begin{array}{ccccc}
X & d & T_d & Y_d & Z_d \\
 & & \\
\G(p,n) & \le\min(p,n-p) & \G(p-d,n-2d) & \G(d,2d) & 
\G(p-d,p)\times\G(n-p-d,n-p) \\
\G_{\omega}(n,2n) & \le n & \G_{\omega}(n-d,2n-2d) & \G_{\omega}(d,2d)
& \G(n-d,n) \\
\G_Q(n,2n) & \le [n/2] & \G_{Q}(n-2d,2n-4d) &\G_{Q}(2d,4d) & \G(n-2d,n)\\
\OO\PP^2 & 1 & \PP^5 & \PP^1 & G_Q(5,10) \\
 & 2 & {\rm pt} & \QQ^8 & \QQ^8\\
E_7/P_7 & 1 & \QQ^{10} & \PP^1 & \OO\PP^2 \\
 & 2 & \PP^1 & \QQ^{10} & \OO\PP^2 \\
 & 3 & {\rm pt} & E_7/P_7 & {\rm pt} 
\end{array}$$ 

\begin{lemm}
\label{intervalles}
For any $d\le d_{max}$, the Bruhat interval $[T_d,Y_d^*]$ in $W_X$
can be identified with the Hasse diagram of a minuscule homogeneous 
variety $Z_d$.  
\end{lemm}

\proof
For any Schubert subvariety $X(w)$ of $X$, we denoted by 
$F_d(\wh)=q_d(p_d^{-1}(X(w))$ the corresponding Schubert subvariety 
of $F_d$. 
The map $w\mapsto\wh$ identifies the Bruhat interval $[T_d,Y_d^*]$
in $W_X$ with the Bruhat interval $[Z_d^*,F_d]$ in $W_{F_d}$. 
By Poincar{\'e} duality on $F_d$ and $Z_d$, this interval is 
canonically identified with $[pt,Z_d]$, which coincides with the
Hasse diagram of $Z_d$.\qed

\medskip
Moreover, each interval $W_d$, being isomorphic with the Hasse 
diagram of the smooth homogeneous variety $Z_d$, is endowed with a natural 
involution induced by Poincar{\'e} duality on the latter. We get 
a global involution of $W_X$ which we denote by $\iota$. In particular
$\iota(T_d)=Y_d^*$, and the fact that $I_d(Y_d^*,\{pt\},T_d)=1$ 
implies that for any $u\in W_X$,
$${\rm codim} (X(u))+ {\rm codim} (X(\iota(u)))=\delta(u)c_1(X).$$   

\section{The quantum product by the class of a point}

In this section we want to compare the intersection in $F_d$ 
of two cells in $[Z_d^*,F_d]$
with an intersection in $Z_d$. Our first result will hold for 
arbitrary homogeneous spaces.

\smallskip
Let $X$ be a homogeneous space
and $Y \subset X$ a Schubert subvariety. Suppose that $Y$ is smooth, 
and even homogeneous. The cell decomposition of $X$ defined by the 
Schubert cells, gives a cell decomposition of $Y$ if we  consider 
only those Schubert cells that are contained in $Y$. In particular
this yields a natural inclusion $W_Y\hookrightarrow W_X$. 
We denote $p_X : W_X \rightarrow W_X$ 
the Poincar{\'e} duality on $X$ and $p_Y : W_Y \rightarrow W_Y$ 
the Poincar{\'e} duality on $Y$. 

Consider the submodule $\Z [Y^*,X]$ of $H^*(X)$ generated by the 
Schubert classes in 
$[Y^*,X]$. This module has a
natural algebra structure: if $[X(u)] \cup [X(v)] = 
\sum_{w \in W_X} c_{u,v}^w [X(w)]$ in $X$, set 
$$[X(u)]\cdot [X(v)] = \sum_{w \in [Y^*,X]} c_{u,v}^w [X(w)].$$
Let $j$ denote the inclusion $Y \rightarrow X$; 
the morphism of modules $j_*:H_*(Y) \rightarrow H_*(X)$ is injective
and the dual  morphism of algebras $j^*:H^*(X) \rightarrow H^*(Y)$ 
is surjective. The  image of $j_*$ is   $\Z [pt,Y]$; 
we denote by  $j_*^{-1}$ the inverse
map $\Z [pt,Y]\ra H_*(Y)$.

\begin{prop}\label{j*}
The restriction of $j^*$ to $\Z[Y^*,X]$ is an algebra isomorphism
with $H^*(Y)$. Explicitely,  it equals $p_Y \circ j_*^{-1} \circ p_X$.
\end{prop}

\proof The first assertion is an immediate consequence of 
\cite[lemma 3.12]{cmp}. To prove the second one, 
let $C \in \Z [pt,Y] \subset H_*(X)$, and let $D \in H^*(X)$.
Since $j_*$ and $j^*$ are adjoint, we have: 
$$
\begin{array}{rcl}
(j_*p_Yj^*p_X(C),D)_X & = & (p_Yj^*p_X(C),j^*D)_Y\\
& = & \delta_{j^*p_X(C),j^*D}\\
& = & \delta_{p_X(C),D} = (C,D)_X
\end{array}
$$
The second and the third equality follow from the fact that
$j^*$ maps bijectively Schubert classes to Schubert classes.
We therefore have proved that $C=j_*p_Yj^*p_X(C)$ for all $C$ 
in $\Z [pt,Y]$ or, equivalently, that $j^*D=p_Yj_*^{-1}p_X(D)$
for all $D$ in $\Z[Y^*,X]$.\qed

\begin{coro}\label{prolisse}
 Let $Y\subset X$ be a homogeneous Schubert subvariety 
of some rational homogeneous manifold. Let $Z$ be any other Schubert
subvariety of $X$. Then 
$$[Y]\cup [Z] = \left\{\begin{array}{ll} p_Xp_Y[Z] & if\;Z\supset Y^*,
    \\ 0 & otherwise.\end{array}\right.$$
\end{coro}

We come back to the (co)minuscule setting.   
Our main examples of homogeneous Schubert subvarieties are our
varieties $T_d$ and $Y_d$ in $X$ and $Z_d$ in $F_d$. Applying the
previous result to the latter varieties, we get the following
statement.

\begin{theo}\label{point}
For any $u\in W_X$, we have
$$\sigma({\rm pt})*\sigma(u)=q^{\delta(u)}\sigma(p\iota (u)).$$
\end{theo}

\proof 
Let $u,v \in W_X$. If $I_d({\rm pt},X(u),X(v)) \not = 0$, then, by
corollary 3.21 in \cite{cmp}, we have $u \leq w_{Y_d^*}$ and
$v \leq w_{Y_d^*}$. Since 
$I_d({\rm pt},X(u),X(v))=I_0(Z_d,F_d(\uh),F_d(\vh)),$
by \cite[lemma 3.12]{cmp}, we also have $\uh,\vh \geq Z_d^*$, which is
equivalent to $u,v \geq T_d$ by the proof of Lemma \ref{intervalles}. 
We conclude that $I_d({\rm pt},X(u),X(v)) \not = 0$ implies 
that $X(u),X(v) \in [T_d,Y_d^*]$,which is equivalent to 
$d=\delta(u)=\delta(v)$.

Denote by 
$j_d : Z_d \hookrightarrow F_d$ the inclusion map. We have
$$I_0(Z_d,F_d(\vh),F_d(\uh))=\int_{Z_d}j_d^*[F_d(\vh)]\cup j_d^*[F_d(\uh)].$$
According to proposition \ref{j*}, 
$\int_{Z_d}j_d^*[F_d(\vh)]\cup j_d^*[F_d(\uh)]$ equals
one if $v=\iota(u)$ and equals zero otherwise. 
The theorem follows.
\qed

\section{Strange duality}

In this section we prove the strange duality property of the quantum
Chow ring stated in Theorem \ref{SD}. For convenience we will in fact
prove a slightly different (but equivalent) statement, asserting the 
existence of a duality mapping $q$ to $q^{-1}$. With the notable 
exception of Grassmannians, this is only possible if we extend scalars 
a little bit. 

We will need a case by case analysis to fix the 
coefficients $\zeta(w)$ in the next statement.

\begin{theo}\label{qinv}
For any (co)minuscule homogeneous variety $X$, one can find an 
algebraic number $\kappa$, and a map $\zeta : W_X\ra
\ZZ[\kappa]$, such that the correspondence 
$$q\mapsto q^{-1}, \qquad 
\sigma(w)\mapsto \zeta(w) q^{-\delta(w)}\sigma(\iota(w))$$
defines a ring involution of $QA^*(X)_{loc}[\kappa]$.
\end{theo}

Clearly, such an involution changes the degree into its opposite, and 
the fact that it is involutive is tantamount to the simple relation 
$\zeta\circ\iota=\zeta^{-1}$. The hard part of the statement is that
it is compatible with the quantum product. 

To get rid of $\kappa$, there remains to 
compose the previous involution with a degree automorphism 
(that multiplies a degree $d$ class by $t^d$ for some $t$). 
It is then a routine but again case by case check that 
we obtain Theorem \ref{SD}. That is, we check that 
$$\zeta(w)=y(w)y(s_{\alpha_0})^{-\frac{\ell(w)}{c_1(X)}}.$$
But we have no convincing 
explanation of why such a formula should hold true. 

\subsection{Quadrics}

We begin with the easy example of quadrics. We will see that 
already in that case fixing the coefficients involves some 
subtelties. 

\subsubsection*{Even dimensions}
Let $\QQ^{2m}$ be a quadric of dimension $2m$, acted on by
$SO_{2m+2}$. The simple root that defines the corresponding maximal
parabolic subgroup is $\alpha_1$. The Hasse diagram is the following:

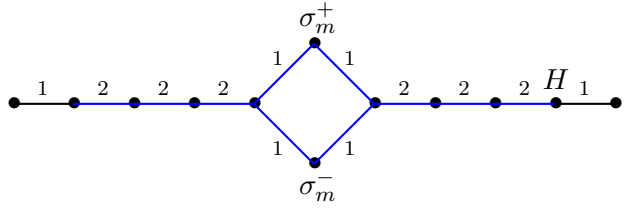
\begin{center}
\psset{unit=4mm}
\psset{xunit=4mm}
\psset{yunit=4mm}
\begin{pspicture*}(20,8)(-10,0)
\multiput(-4.2,3.8)(2,0){5}{$\bullet$}
\multiput(7.8,3.8)(2,0){5}{$\bullet$}
\multiput(5.8,1.8)(0,4){2}{$\bullet$}
\psline(-4,4)(-2,4)
\psline[linecolor=blue](-2,4)(4,4)
\psline[linecolor=blue](4,4)(6,6)
\psline[linecolor=blue](4,4)(6,2)
\psline[linecolor=blue](6,6)(8,4)
\psline[linecolor=blue](6,2)(8,4)
\psline[linecolor=blue](8,4)(14,4)
\psline(14,4)(16,4)
\put(5.5,1){$\sigma_m^-$}\put(5.5,6.6){$\sigma_m^+$}
\put(13.6,4.5){$H$}
\put(-3.2,4.3){$\scriptstyle{1}$}
\put(7,5.3){$\scriptstyle{1}$}
\put(4.6,5.3){$\scriptstyle{1}$}
\put(7,2.3){$\scriptstyle{1}$}
\put(4.6,2.3){$\scriptstyle{1}$}
\multiput(-1.2,4.3)(2,0){3}{$\scriptstyle{2}$}
\multiput(8.8,4.3)(2,0){3}{$\scriptstyle{2}$}
\put(14.8,4.3){$\scriptstyle{1}$}
\end{pspicture*}
\end{center}

\noindent 
Indeed, as is well known, there is one Schubert class $\sigma_k$ in each 
degre $k\ne m$, and two Schubert classes $\sigma_m^\pm$ in middle 
degree, defined by the two families of maximal linear spaces in
$\QQ^{2m}$. Of course $\sigma_1=H$ is the hyperplane class. 
For the classical intersection product, we have 
$H^k=\sigma_k$ for $k<m$, $H^m=\sigma_m^++\sigma_m^-$, and 
$H^k=2\sigma_k$ for $k>m$. Since $c_1(\QQ^{2m})=2m$, these formulas 
remain valid in the quantum Chow ring, except in maximal degree. 
In fact the quantum Chevalley formula gives the two identities
$$\sigma_{2m-1}*H=\sigma_{2m}+q, \qquad \sigma_{2m}*H=qH.$$
Note that $\sigma_{2m-1}$ is the class of a line and coincides with
$[T_1]$. Also $\sigma_{2m}$ is the class $[pt]$ of a point. 

\begin{prop} The quantum Chow ring $QA^*(\QQ^{2m})$ is determined by the 
formulas 
$$\sigma_m^+*\sigma_m^-=[pt], \qquad \sigma_m^+*\sigma_m^+=
\sigma_m^-*\sigma_m^-=q.$$
\end{prop}

Note that these relations imply that $H^{2m}=2[pt]+2q$ and 
$H^{2m+1}=4qH$. Now it is easy to check that 
$QA^*(\QQ^{2m})_{loc}$ has a ring involution given by 
$$q\mapsto 1/16q, \qquad H^k\mapsto H^{2m-k}/4q, 
\qquad [pt]\mapsto [pt]/16q^2,$$
(where the central formula holds for $1\le k\le 2m-1$), and in
middle degree  by $\sigma_m^\pm\mapsto\sigma_m^\mp/4q$.
Multiplying each degree $d$ class by $2^d$, we get a slightly
different involution given by 
$$\begin{array}{rl}
 q &\mapsto 2^{2m-4}/q, \\
 \sigma_k &\mapsto 
 2^{k-1}\sigma_{2m-k}/q \quad \mathrm{for}\; 0<k<m, \\
 \sigma_m^\pm &\mapsto 2^{m-2}\sigma_m^\mp, \\
 \sigma_k &\mapsto 2^{k-3}\sigma_{2m-k}/q \quad \mathrm{for}\;m<k<2m, \\
 \hspace{1cm}[pt] &\mapsto 2^{2m-4}[pt]/q^2.
\end{array}$$
This is precisely the statement of Theorem \ref{SD} for $\QQ^{2m}$,
as one can readily check from the Hasse diagram above. Indeed, we have
labeled the edges by the coefficients whose products, taken from the 
rightmost end of the diagram, give the coefficients $y(w)$. 
The corresponding coefficients for Theorem \ref{qinv} are
the following, where $\zeta(m)$ accounts for the two mid-dimensional 
classes:  
$$\zeta(k)=4^{\frac{k-m}{m}}\quad \mathrm{for}\;0<k<2m, \qquad 
\zeta(0)=\zeta(2m)=1,$$

\subsubsection*{Odd dimensions}

An odd dimensional quadric $\QQ^{2m-1}$ is cominuscule but not
minuscule. It has exactly one Schubert class $\sigma_k$ in each 
dimension $k<2m$. Again   $\sigma_1=H$ is the hyperplane class. 
The Hasse diagram is the simplest possible one; as in the even 
dimensional case we have labeled the edges by the coefficients 
whose products, taken from the rightmost end,
 give the $y(w)$ of Theorem \ref{SD}. 

\begin{center}
\psset{unit=4mm}
\psset{xunit=4mm}
\psset{yunit=4mm}
\begin{pspicture*}(20,6)(-10,3)
\multiput(-4.2,3.8)(2,0){10}{$\bullet$}
\psline(-4,4)(-2,4)
\psline[linecolor=blue](-2,4)(12,4)
\psline(12,4)(14,4)
\put(11.6,4.5){$H$}
\put(-3.2,4.3){$\scriptstyle{1}$}
\multiput(-1.2,4.3)(2,0){3}{$\scriptstyle{2}$}
\put(4.8,4.5){$\scriptstyle{\frac{1}{2}}$}
\multiput(6.8,4.3)(2,0){3}{$\scriptstyle{2}$}
\put(12.8,4.3){$\scriptstyle{1}$}
\end{pspicture*}
\end{center}

For the classical intersection product, we have 
$H^k=\sigma_k$ for $k<m$, 
$H^k=2\sigma_k$ for $k\ge m$. Since $c_1(\QQ^{2m-1})=2m-1$, these 
formulas hold in the quantum Chow ring, except in maximal degree. 
In fact the quantum Chevalley formula gives the two identities
$$\sigma_{2m-2}*H=\sigma_{2m-1}+q, \qquad \sigma_{2m-1}*H=qH.$$
This is enough to determine the quantum product. 
The class $\sigma_{2m-2}$ is that of a line and coincides with
$[T_1]$. Also $\sigma_{2m-1}$ is the class $[pt]$ of a point. 

Again it is easy to check that 
$QA^*(\QQ^{2m-1})_{loc}$ has a ring involution given by almost 
the same formulas that in even dimensions,  
$$q\mapsto 1/16q, \qquad H^k\mapsto H^{2m-k-1}/4q, 
\qquad [pt]\mapsto [pt]/16q^2,$$
(where the central formula holds for $1\le k\le 2m-2$).
Multiplying each degree $d$ class by $2^d$, we get the slightly
different involution given by 
$$\begin{array}{rl}
 q &\mapsto 2^{2m-5}/q, \\
 \sigma_k &\mapsto 
 2^{k-1}\sigma_{2m-k-1}/q \quad \mathrm{for}\; 0<k<m, \\
 \sigma_{k} &\mapsto 2^{k-3}\sigma_{2m-k-1}/q \quad \mathrm{for}\;m\le k<2m-1, \\
 \hspace{1cm}[pt] &\mapsto 2^{2m-5}[pt]/q^2.
\end{array}$$
This is in perfect agreement with Theorem \ref{SD}. 

\subsection{Grassmannians}

Now suppose $X=\G(p,n)$ is a Grassmannian. Then $W_X$ identifies with the 
set of partitions inscribed in the rectangle $p\times (n-p)$, and 
$d_{max}=\min(p,n-p)$.  For any $d\subset d_{max}$, $Y_d$ is a
Grassmannian $\G(d,2d)$, and $T_d$ is another Grassmannian
$\G(p-d,n-2d)$. In particular the partition $w_{Y_d^*}$ is just a
square of size $d$, while the partition  $w_{T_d}$ is the complement 
of a rectangle of size $(p-d)\times (n-p-d)$. 

We deduce that for any partition $\lambda\in W_X$, the degree
$d=\delta(\lambda)$ is the size of the biggest square contained 
in $\lambda$. Clearly the interval $[T_d,Y_d^*]$ in $W_X$
identifies with the Hasse diagram of the product $\G(d,p-d)
\times \G(d,n-p-d)$. The partition $\iota(\lambda)$ is deduced from
$\lambda$ by taking the complementary partitions in the SW and NE 
rectangles.

In this case Theorem \ref{qinv} appears in \cite{Po}, Theorem 7.5.
See also \cite{H} where the involution is interpreted in terms of 
complex conjugation.

\subsection{Lagrangian Grassmannians}

Let $X=\G_{\omega}(n,2n)$ be a Lagrangian Grassmannian. Recall that 
$W_X$ identifies with the set of strict partitions
$\lambda\subset\rho_n$. In order to simplify notations we let 
$\sigma_{\lambda}$ denote the class of the Schubert subvariety 
Poincar{\'e} dual to $X(\lambda)$. Its degree is the sum $|\lambda|$
of the parts of $\lambda$. Morevoer, $d=\delta(\lambda)$ is simply the 
number of (non zero) parts of $\lambda$, usually called the 
length and denoted $\ell(\lambda)$. In particular $d_{max}=n$. 

Mapping $\lambda$ to
$(\lambda_1-d,\ldots,\lambda_d-1)$ we get a 
partition inscribed in the rectangle $d\times (n-d)$. This identifies
the interval $W_d$ of $W_X$ with the Hasse diagram of the Grassmannian
$\G(d,n)$. Appying Poincar{\'e} duality for that Grassmannian we dedude
that the involution $\iota$ is given by 
$$\iota (\lambda_1,\ldots, \lambda_{\ell(\lambda)})=(
n+1-\lambda_{\ell(\lambda)},\ldots, n+1-\lambda_1).$$

\begin{example} For $n=5$ we have drawn below the Hasse diagram 
of $\G_{\omega}(5,10)$ and its partition into six disjoint
intervals. More precisely we have drawn in black the arrows between
different intervals, and in the same color the arrows inside a given
interval. Those are the Hasse diagrams of a point, of $\PP^4$ and 
of the Grassmannian $\G(2,5)$, all appearing twice and
symmetrically.

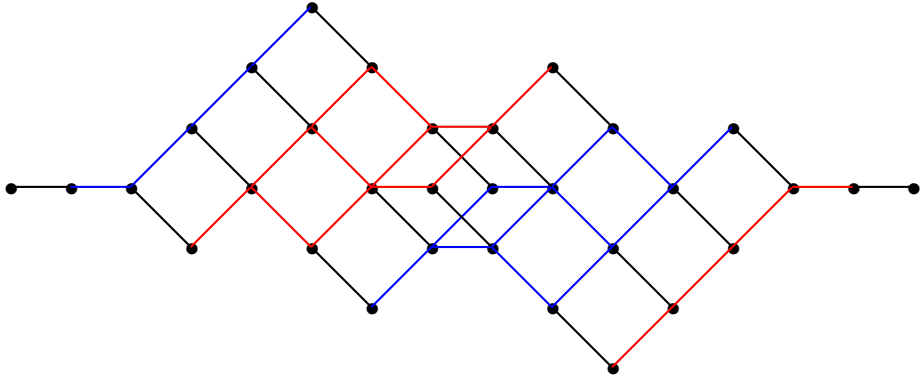
\begin{center}
\psset{unit=4mm}
\psset{xunit=4mm}
\psset{yunit=4mm}
\begin{pspicture*}(40,15)(-10,0)
\multiput(-4.2,7.7)(2,0){3}{$\bullet$}
\multiput(1.8,9.7)(4,0){3}{$\bullet$}
\multiput(11.8,9.7)(4,0){3}{$\bullet$}
\multiput(1.8,5.7)(4,0){3}{$\bullet$}
\multiput(11.8,5.7)(4,0){3}{$\bullet$}
\put(7.8,3.7){$\bullet$}
\put(15.8,1.7){$\bullet$}
\multiput(13.8,3.7)(4,0){2}{$\bullet$}
\multiput(3.8,11.7)(4,0){2}{$\bullet$}
\multiput(3.8,7.7)(4,0){3}{$\bullet$}
\multiput(9.8,7.7)(4,0){3}{$\bullet$}
\put(13.8,11.7){$\bullet$}
\multiput(21.8,7.7)(2,0){3}{$\bullet$}
\put(5.8,13.7){$\bullet$}
\psline(-4,8)(-2,8)
\psline[linecolor=blue](-2,8)(0,8)
\psline(0,8)(2,6)
\psline[linecolor=blue](0,8)(6,14)
\psline[linecolor=red](2,6)(8,12)
\psline(2,10)(4,8)\psline[linecolor=red](4,8)(6,6)\psline(6,6)(8,4)
\psline(4,12)(6,10)\psline[linecolor=red](6,10)(8,8)\psline(8,8)(10,6)
\psline[linecolor=red](6,6)(10,10)
\psline(6,14)(8,12)\psline[linecolor=red](8,12)(10,10)\psline(10,10)(12,8)
\psline[linecolor=blue](8,4)(12,8)
\psline[linecolor=red](8,8)(10,8)
\psline[linecolor=blue](10,6)(12,6)
\psline[linecolor=red](10,10)(12,10)
\psline[linecolor=blue](12,8)(14,8)
\psline(10,8)(12,6)\psline[linecolor=blue](12,6)(14,4)\psline(14,4)(16,2)
\psline[linecolor=red](10,8)(14,12)
\psline[linecolor=blue](12,6)(16,10)
\psline(12,10)(14,8)\psline[linecolor=blue](14,8)(16,6)\psline(16,6)(18,4)
\psline(14,12)(16,10)\psline[linecolor=blue](16,10)(18,8)\psline(18,8)(20,6)
\psline[linecolor=blue](14,4)(20,10)
\psline[linecolor=red](16,2)(22,8)
\psline(20,10)(22,8)\psline(24,8)(26,8)
\psline[linecolor=red](22,8)(24,8)
\end{pspicture*}
\end{center}
\end{example}

\begin{prop} 
For $\lambda\subset\rho_n$ a strict partition, let 
$z(\lambda):=\ell(\lambda)-\frac{2|\lambda|}{n+1}$. 
Then $z(\iota(\lambda))=z(\lambda^*)=-z(\lambda)$, and Theorem
\ref{qinv} holds with 
$$\kappa=2^{-\frac{2}{n+1}} \qquad and \qquad
\zeta(\lambda)=2^{z(\lambda)}.$$
\end{prop}

\proof We just need to prove that our involution is compatible 
with the quantum Pieri formula for $\sigma_{\lambda}*\sigma_k$, for
each $k$. We have $\iota(\sigma_k)=q^{-1}\sigma_{n+1-k}$, so we need to 
compare $\sigma_{\lambda}*\sigma_k$ with $\sigma_{ \iota(\lambda)}*
\sigma_{n+1-k}$.

We use the quantum version of Pieri's rule as stated in \cite{BKT},
Theorem 3:
$$\sigma_{\lambda}*\sigma_k=\sum_{\mu}2^{N(\lambda,\mu)}\sigma_{\mu}
+q\sum_{\nu}2^{N'(\nu,\lambda)}\sigma_{\nu},$$
where the first sum is over all strict partitions $\mu\supset\lambda$
with $|\mu|=|\lambda|+k$, such that the complement $\mu/\lambda$ is a
horizontal strip, and the second sum is over all strict partitions 
$\nu\subset\lambda$ with $|\nu|=|\lambda|-(n+1-k)$, such that the 
$\lambda/\nu$ is a horizontal strip. Moreover $N(\lambda,\mu)$ denotes
the number of connected components of $\mu/\lambda$ which do no meet
the first column, and $N'(\nu,\lambda)$ is one less than the total 
number of connected components of $\lambda/\nu$. (By definition, two 
boxes are connected if they share an edge or a vertex.) Applying our
involution, we get 
\begin{eqnarray}\label{1} \nonumber
\iota(\sigma_{\lambda}*\sigma_k) = 
&q^{-\ell(\lambda)}\sum_{\ell(\mu)=\ell(\lambda)}2^{z(\mu)+N(\lambda,\mu)}
\sigma_{\iota(\mu)}
+q^{-\ell(\lambda)-1}\sum_{\ell(\mu)=\ell(\lambda)+1}
2^{z(\mu)+N(\lambda,\mu)}\sigma_{\iota(\mu)} \\
  &+q^{-\ell(\lambda)-1}\sum_{\ell(\nu)=\ell(\lambda)}
2^{z(\nu)+N'(\nu,\lambda)}\sigma_{\iota(\nu)}
+q^{-\ell(\lambda)}\sum_{\ell(\nu)=\ell(\lambda)-1}
2^{z(\nu)+N'(\nu,\lambda)}\sigma_{\iota(\nu)}.
\end{eqnarray}
Similarly, we deduce again from Pieri's rule that 
\begin{eqnarray}\label{2} \nonumber
\iota(\sigma_{\lambda})*\iota(\sigma_k) = &2^{z(\lambda)+z(k)}
\Big(q^{-\ell(\lambda)-1}\sum_{\ell(\alpha)=\ell(\lambda)}
2^{N(\iota(\lambda),\alpha)}\sigma_{\alpha}
+q^{-\ell(\lambda)-1}\sum_{\ell(\alpha)=\ell(\lambda)+1} 
 2^{N(\iota(\lambda),\alpha)}\sigma_{\alpha} \\
 &+q^{-\ell(\lambda)}\sum_{\ell(\beta)=\ell(\lambda)}
2^{N'(\beta,\iota(\lambda))}\sigma_{\beta}
+q^{-\ell(\lambda)}\sum_{\ell(\beta)=\ell(\lambda)-1}
2^{N'(\beta,\iota(\lambda))}\sigma_{\beta}\Big).
\end{eqnarray}
We claim that the four partial sums on the right hand sides of the 
two identities above correspond pairwise. 

Consider the first term on the right hand side of (\ref{2}). 
Here $\alpha$ is a strict partition containing $\iota(\lambda)$, of size 
$|\alpha|=|\iota(\lambda)|+n+1-k$, such that $\alpha/\iota(\lambda)$
is a horizontal strip. If we let $\nu=\iota(\alpha)$, we get that 
$\nu$ is contained in $\lambda$, $\ell(\nu)=\ell(\alpha)$ and
$|\nu|=|\lambda|-(n+1-k)$, and $\lambda/\nu\simeq
\alpha/\iota(\lambda)$ is a horizontal strip. Moreover, the fact that 
$\ell(\alpha)=\ell(\lambda)$ means that $\alpha/\iota(\lambda)$ does
not meet the first column, so
$N(\iota(\lambda),\alpha)=N'(\nu,\lambda)+1$. On the other hand, 
$$z(\nu)=\ell(\lambda)-\frac{2}{n+1}(|\lambda|-(n+1-k))=z(\lambda)+z(k)+1,$$
and therefore
$z(\nu)+N'(\nu,\lambda)=z(\lambda)+z(k)+N(\iota(\lambda),\alpha)$. We
conclude that the first partial sum of (\ref{2}) coincides exactly 
with the third partial sum of (\ref{1}).

For the second term on the right hand side of (\ref{2}), 
the difference with the first term is that 
$\ell(\alpha)=\ell(\lambda)+1$, which means that
$\alpha/\iota(\lambda)$ does meet the first column. Thus
$N(\iota(\lambda),\alpha)=N'(\nu,\lambda)$. On the other hand,
$$z(\nu)=\ell(\lambda)+1-\frac{2}{n+1}(|\lambda|+(n+1)-(n+1-k))
=z(\lambda)+z(k),$$
and therefore we get again 
$z(\nu)+N'(\nu,\lambda)=z(\lambda)+z(k)+N(\iota(\lambda),\alpha)$.
We conclude that the second partial sum of (\ref{2}) coincides exactly 
with the second partial sum of (\ref{1}). 

Now consider the third term on the right hand side of (\ref{2}).  
Here $\beta$ is a strict partition contained in
$\iota(\lambda)$, with $\ell(\beta)=\ell(\lambda)$ and $|\beta|=
|\iota(\lambda)|-k$, such that $\iota(\lambda)/\beta$ is a 
horizontal strip. If we let $\beta=\iota(\mu)$, we get 
$\iota(\lambda)/\beta\simeq\mu/\lambda$, so $|\mu|=|\lambda|+k$ and 
$N'(\beta,\iota(\lambda))=N'(\lambda,\mu)=N(\lambda,\mu)-1$. Again
we deduce that  $z(\lambda)+z(k)+N'(\beta,\iota(\lambda))=z(\mu)+
N(\lambda,\mu)$,  so that the third partial sum of (\ref{2})
coincides exactly with the first partial sum of (\ref{1}). 

Finally, consider the third term on the right hand side of (\ref{2}). 
Here $\beta$ is as before except that
$\ell(\beta)=\ell(\lambda)-1$, so that if we let $\nu=\iota(\beta)$,
then  $N'(\beta,\iota(\lambda))=N'(\nu,\lambda)$. Again we conclude 
that the fourth partial sum of (\ref{2}) coincides exactly with the 
fourth partial sum of (\ref{1}). 
 
This concludes the proof. \qed

\medskip As we have already mentionned, it is possible to get rid of
the root of two by composing with a degree automorphism. We get:

\begin{theo}\label{inv-LG}
The correspondence 
$$q\mapsto \frac{4}{q}, \hspace{15mm}
\sigma_{\lambda}\mapsto \Big(\frac{2}{q} \Big)^{\ell(\lambda)}
\sigma_{\iota(\lambda)} \quad for\;\lambda\subset\rho_n,$$
defines a ring involution of $QA^*(G_{\omega}(n,2n))_{loc}$.
\end{theo}

\subsection{Orthogonal Grassmannians}

Let $X=\G_Q(n+1,2n+2)$ be a Lagrangian Grassmannian. Recall that 
$W_X$ identifies again with the set of strict partitions
$\lambda\subset\rho_n$. Again we denote by $\sigma_{\lambda}$ the
Schubert  class Poincar{\'e} dual to $[X(\lambda)]$. Its degree 
is $|\lambda|$, and $\delta(\lambda)=d$ if 
the number of (non zero) parts of $\lambda$ is $2d$ or $2d-1$. 
In particular $d_{max}=[n/2]$. 

Mapping $\lambda$ to
$(\lambda_1-2d+1,\ldots,\lambda_{2d})$ (where we let $\lambda_{2d}=0$
if $\ell(\lambda)=2d-1$), we get a 
partition inscribed in the rectangle $2d\times (n-2d+1)$. This identifies
the interval $W_d$ of $W_X$ with the Hasse diagram of the Grassmannian
$\G(2d,n+1)$. Appying Poincar{\'e} duality for that Grassmannian we dedude
that the involution $\iota$ is given by 
$$\iota (\lambda_1,\ldots, \lambda_{2\delta(\lambda)})=(
n-\lambda_{2\delta(\lambda)},\ldots, n-\lambda_1).$$

\begin{example} For $n=5$ we have drawn below the Hasse diagram 
of $\G_Q(6,12)$ and its partition into four disjoint
intervals. Note that the diagram is the same as for
$\G_{\omega}(5,10)$ but the partition is different. 
Indeed we have only four intervals in that case, 
isomorphic with the Hasse diagrams of a point and of $G(2,6)$
appearing twice symmetrically. (Beware that this symmetry is specific
to the case where $n$ is odd.)
 
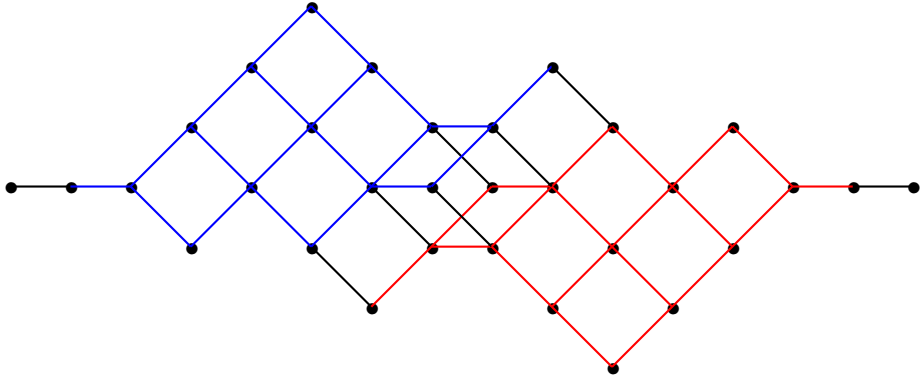
\begin{center}
\psset{unit=4mm}
\psset{xunit=4mm}
\psset{yunit=4mm}
\begin{pspicture*}(40,15)(-10,0)
\multiput(-4.2,7.7)(2,0){3}{$\bullet$}
\multiput(1.8,9.7)(4,0){3}{$\bullet$}
\multiput(11.8,9.7)(4,0){3}{$\bullet$}
\multiput(1.8,5.7)(4,0){3}{$\bullet$}
\multiput(11.8,5.7)(4,0){3}{$\bullet$}
\put(7.8,3.7){$\bullet$}
\put(15.8,1.7){$\bullet$}
\multiput(13.8,3.7)(4,0){2}{$\bullet$}
\multiput(3.8,11.7)(4,0){2}{$\bullet$}
\multiput(3.8,7.7)(4,0){3}{$\bullet$}
\multiput(9.8,7.7)(4,0){3}{$\bullet$}
\put(13.8,11.7){$\bullet$}
\multiput(21.8,7.7)(2,0){3}{$\bullet$}
\put(5.8,13.7){$\bullet$}
\psline(-4,8)(-2,8)
\psline[linecolor=blue](-2,8)(0,8)
\psline[linecolor=blue](0,8)(2,6)
\psline[linecolor=blue](0,8)(6,14)
\psline[linecolor=blue](2,6)(8,12)
\psline[linecolor=blue](2,10)(6,6)\psline(6,6)(8,4)
\psline[linecolor=blue](4,12)(8,8)\psline(8,8)(10,6)
\psline[linecolor=blue](6,6)(10,10)
\psline[linecolor=blue](6,14)(10,10)\psline(10,10)(12,8)
\psline[linecolor=red](8,4)(12,8)
\psline[linecolor=blue](8,8)(10,8)
\psline[linecolor=red](10,6)(12,6)
\psline[linecolor=blue](10,10)(12,10)
\psline[linecolor=red](12,8)(14,8)
\psline(10,8)(12,6)\psline[linecolor=red](12,6)(16,2)
\psline[linecolor=blue](10,8)(14,12)
\psline[linecolor=red](12,6)(16,10)
\psline(12,10)(14,8)\psline[linecolor=red](14,8)(18,4)
\psline(14,12)(16,10)\psline[linecolor=red](16,10)(20,6)
\psline[linecolor=red](14,4)(20,10)
\psline[linecolor=red](16,2)(22,8)
\psline[linecolor=red](20,10)(22,8)\psline(24,8)(26,8)
\psline[linecolor=red](22,8)(24,8)
\end{pspicture*}
\end{center}
\end{example}

\begin{prop} 
For $\lambda\subset\rho_n$ a strict partition, let 
$z(\lambda):=\frac{2|\lambda|}{n}-(\ell(\lambda)+\delta_{\lambda_1,n})$.
Then $z(\iota(\lambda))=z(\lambda^*)=-z(\lambda)$, and Theorem
\ref{qinv} holds with 
$$\kappa=2^{-\frac{2}{n}} \qquad and \qquad
\zeta(\lambda)=2^{z(\lambda)}.$$
\end{prop}

\proof The proof is notably different from that of Theorem
\ref{inv-LG}, since we have 
$$\iota(\sigma_k)=q^{-1}\sigma_{n,n-k}=q^{-1}\sigma_n*\sigma_{n-k}.$$
Again we check that our formula is compatible with the quantum 
version of Pieri's rule as stated in \cite{BKT},
Theorem 6:
$$\sigma_{\lambda}*\sigma_k=\sum_{\mu}2^{N'(\lambda,\mu)}\sigma_{\mu}
+q\sum_{\nu}2^{N'(\lambda,\nu)}\sigma_{\bar{\nu}},$$
where the first sum is over all strict partitions $\mu\supset\lambda$
with $|\mu|=|\lambda|+k$, such that $\mu/\lambda$ is a
horizontal strip, and the second sum is over all partitions 
$\nu=(n,n,\bar{\nu})\supset\lambda$, with $\bar{\nu}$ strict and 
$|\nu|=|\lambda|+k$, such that  
$\nu/\lambda$ is a horizontal strip. Note that the quantum correction 
is non trivial zero when $\lambda_1=n$. 

We distinguish several cases. 

\medskip\noindent {\sc First case}: $\lambda_1<n$. Then 
\begin{eqnarray}\label{3}
\sigma_{\lambda}*\sigma_k=
\sum_{\ell(\mu)=\ell(\lambda)}2^{N'(\lambda,\mu)}\sigma_{\mu}
+\sum_{\ell(\mu)=\ell(\lambda)+1}2^{N'(\lambda,\mu)}\sigma_{\mu}
\end{eqnarray}
has no quantum correction. 

\smallskip\noindent {\it First subcase}: $\iota(\lambda)_1<n$. This
means that $\ell(\lambda)=2\delta(\lambda)$, so that on the 
right hand side of  (\ref{3}) 
we have $\delta(\mu)=\delta(\lambda)$ in the first 
partial sum and $\delta(\mu)=\delta(\lambda)+1$ in the second one. 
Hence 
\begin{eqnarray}\label{5}
\iota(\sigma_{\lambda}*\sigma_k)=
q^{-\delta(\lambda)}\sum_{\ell(\mu)=\ell(\lambda)}2^{z(\mu)+N'(\lambda,\mu)}
\sigma_{\iota(\mu)}
+q^{-\delta(\lambda)-1}\sum_{\ell(\mu)=\ell(\lambda)+1}2^{z(\mu)+N'(\lambda,\mu)}
\sigma_{\iota(\mu)}.
\end{eqnarray}
On the other hand
$\sigma_{\iota(\lambda)}*\sigma_n=\sigma_{n,\iota(\lambda)}$
has no quantum correction, and the quantum Pieri rule gives 
\begin{eqnarray}\label{6} 
\iota(\sigma_{\lambda})*\iota(\sigma_k) = 2^{z(\lambda)+z(k)}
\Big(q^{-\delta(\lambda)-1}\sum_{\alpha}
2^{N'(\iota(\lambda),\alpha)}\sigma_{n,\alpha}
+q^{-\delta(\lambda)}\sum_{\beta}
2^{N'((n,\iota(\lambda)),\beta)}\sigma_{\bar{\beta}}\Big).
\end{eqnarray}
Consider some $\mu$ in the first sum on the right hand side of 
(\ref{5}), and let $\beta=\iota(\mu)$. We claim that 
$$N'((n,\iota(\lambda)),\beta)=N'(\lambda,\mu)+1-\delta_{\mu_1,n}.$$
The following picture should help to see this. We have represented 
the partition $\lambda$ in thick lines, so that $\iota(\lambda)$ is
its complement (reversed) in the rectangle $\ell(\lambda)\times n$. 
We have added a line a the bottom of this rectangle to represent 
$(n,\iota(\lambda))$ (again reversed). The $\bullet$'s represent
$\beta/(n,\iota(\lambda))$ (reversed), a horizontal strip  inside $\lambda$ --
except possibly if $\ell(\beta)=\ell(\lambda)+2$, in which case there 
are some $\bullet$'s on the line above the first line of $\lambda$. 
The $\circ$'s represent $\mu/\lambda$, again a horizontal strip. On
each line they complement the $\bullet$'s of the line above. So we
start with a connected component of $\bullet$'s on the SW corner of
the picture, and going NE we  successively meet the connected components
of $\circ$'s and $\bullet$'s. We have thus the same number of
components, or one more for the $\bullet$'s if we end by one of 
these at the NE corner. This is the case if and only if $\mu_1<n$,
so our claim follows. 

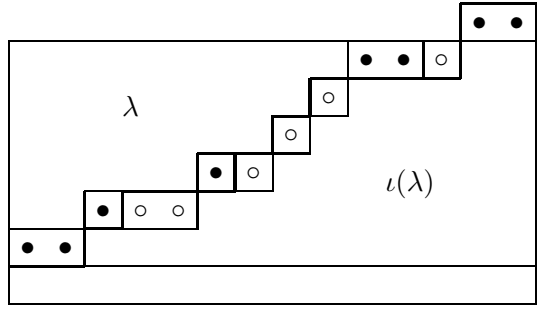
\begin{center}
\setlength{\unitlength}{5mm}
\begin{picture}(20,8)(-3,0)
\put(0,0){\line(1,0){14}}
\put(0,1){\line(1,0){14}}
\put(0,7){\line(1,0){14}}
\put(0,0){\line(0,1){7}}
\put(14,0){\line(0,1){7}}
\put(0,2){\line(1,0){2}}
\put(2,2){\line(0,1){1}}
\put(2,3){\line(1,0){1}}
\put(5,3){\line(0,1){1}}
\put(5,4){\line(1,0){1}}
\put(9,6){\line(0,1){1}}
\put(12,7){\line(0,1){1}}
\put(14,7){\line(0,1){1}}
\put(12,8){\line(1,0){2}}
\put(.3,1.3){$\bullet$}
\put(1.3,1.3){$\bullet$}
\put(2.3,2.3){$\bullet$}
\put(5.3,3.3){$\bullet$}
\put(9.3,6.3){$\bullet$}
\put(10.3,6.3){$\bullet$}
\put(12.3,7.3){$\bullet$}
\put(13.3,7.3){$\bullet$}
\put(3,2){\line(1,0){2}}
\put(5,2){\line(0,1){1}}
\put(6,3){\line(1,0){1}}
\put(7,3){\line(0,1){1}}
\put(7,4){\line(1,0){1}}
\put(8,4){\line(0,1){1}}
\put(8,5){\line(1,0){1}}
\put(9,5){\line(0,1){1}}
\put(11,6){\line(1,0){1}}
\put(12,6){\line(0,1){1}}
\put(3.3,2.3){$\circ$}
\put(4.3,2.3){$\circ$}
\put(6.3,3.3){$\circ$}
\put(7.3,4.3){$\circ$}
\put(8.3,5.3){$\circ$}
\put(11.3,6.3){$\circ$}
\thicklines
\put(0,1){\line(1,0){2}}
\put(2,1){\line(0,1){1}}
\put(2,2){\line(1,0){1}}
\put(3,2){\line(0,1){1}}
\put(3,3){\line(1,0){3}}
\put(6,3){\line(0,1){1}}
\put(6,4){\line(1,0){1}}
\put(7,4){\line(0,1){1}}
\put(7,5){\line(1,0){1}}
\put(8,5){\line(0,1){1}}
\put(8,6){\line(1,0){3}}
\put(11,6){\line(0,1){1}}
\put(3,5){$\lambda$}
\put(10,3){$\iota(\lambda)$}
\end{picture}
\end{center}

But $z(\lambda)+z(k)-z(\mu)=\delta_{\mu_1,n}-1$, so 
$z(\lambda)+z(k)+N'((n,\iota(\lambda)),\beta)=z(\mu)+N'(\lambda,\mu)$.
We conclude that the first sum of (\ref{5}) coincides with the 
second sum of (\ref{6}). 

Now consider some $\mu$ in the second sum on the right hand side of 
(\ref{5}). Since $\ell(\mu)=2\delta(\lambda)+1$ is odd, we get 
$\iota(\mu)=(n,\alpha)$ for some strict partition $\alpha$. We claim
that $$N'(\iota(\lambda),\alpha)=N'(\lambda,\mu)-\delta_{\mu_1,n}.$$
Again this implies that the second sum of (\ref{5}) coincides with the 
first sum of (\ref{6}), and we are done.

\smallskip\noindent {\it Second subcase}: $\iota(\lambda)_1=n$. Then 
$\ell(\lambda)=2\delta(\lambda)-1$, and in  
(\ref{3}) we always have $\delta(\mu)=\delta(\lambda)$. Hence 
\begin{eqnarray}\label{4}
\iota(\sigma_{\lambda}*\sigma_k)=
q^{-\delta(\lambda)}\sum_{\mu}2^{z(\mu)+N'(\lambda,\mu)}
\sigma_{\iota(\mu)}.
\end{eqnarray}
On the other hand,
$\sigma_{\iota(\lambda)}*\sigma_n=q\sigma_{\iota(\lambda)/n}$, where
the first part of $\iota(\lambda)/n$ is smaller than $n$. In
particular, the quantum Pieri rule for $\sigma_{\iota(\lambda)/n}*
\sigma_{n-k}$ has no quantum correction. We get
\begin{eqnarray} 
\iota(\sigma_{\lambda})*\iota(\sigma_k) = 2^{z(\lambda)+z(k)}
q^{-\delta(\lambda)}\sum_{\alpha}
2^{N'(\iota(\lambda)/n,\alpha)}\sigma_{\alpha},
\end{eqnarray}
the sum being taken over all strict partitions
$\alpha\supset\iota(\lambda)/n$, 
with $|\alpha|=|\iota(\lambda)/n|+n-k=|\iota(\lambda)|-k$, such that 
$\alpha/(\iota(\lambda)/n)$  is a horizontal strip. Let
$\mu=\iota(\alpha)$, and note that
$\delta(\mu)=\delta(\alpha)=\delta(\alpha)$
since $\ell(\iota(\lambda)/n)=\ell(\lambda)$ is odd. 
We claim that 
$$N'(\iota(\lambda)/n,\alpha)=N'(\lambda,\mu)+
\delta_{\ell(\mu),\ell(\lambda)}-\delta_{\mu_1,n}.$$
Since we have
$z(\mu)-z(\lambda)-z(k)=1+\ell(\lambda)-\ell(\mu)-\delta_{\mu_1,n}=
\delta_{\ell(\mu),\ell(\lambda)}-\delta_{\mu_1,n}$, 
we conclude that the coefficients
of $\sigma_{\mu}$ in (\ref{3}) and (\ref{4}) are equal, which is what
we wanted to prove. 

\medskip\noindent {\sc Second case}: $\lambda_1=n$. Then we must take
care of the quantum correction in Pieri's rule. We write 
\begin{eqnarray}\label{7}
\sigma_{\lambda}*\sigma_k=
\sum_{\ell(\mu)=\ell(\lambda)}2^{N'(\lambda,\mu)}\sigma_{\mu}
+\sum_{\ell(\mu)=\ell(\lambda)+1}2^{N'(\lambda,\mu)}\sigma_{\mu}
+q\sum_{\ell(\nu)=\ell(\lambda)}2^{N'(\lambda,\nu)}\sigma_{\bar{\nu}}
+q\sum_{\ell(\nu)=\ell(\lambda)+1}2^{N'(\lambda,\nu)}\sigma_{\bar{\nu}}.
\end{eqnarray}

\smallskip\noindent {\it First subcase}: $\iota(\lambda)_1<n$. Then 
$\delta(\mu)=\delta(\lambda)$ in the first sum of (\ref{7}) but
$\delta(\mu)=\delta(\lambda)+1$ in the second sum, while 
$\delta(\bar{\nu})=\delta(\lambda)-1$ in the third sum and
$\delta(\bar{\nu})=\delta(\lambda)$ in the last one. Thus
\begin{eqnarray}\label{8}\nonumber
\iota(\sigma_{\lambda}*\sigma_k)=
&q^{-\delta(\lambda)}\sum_{\ell(\mu)=\ell(\lambda)}2^{z(\mu)+N'(\lambda,\mu)}
\sigma_{\iota(\mu)}
+q^{-\delta(\lambda)-1}\sum_{\ell(\mu)=\ell(\lambda)+1}
2^{z(\mu)+N'(\lambda,\mu)}\sigma_{\iota(\mu)} \\
 &+q^{-\delta(\lambda)}\sum_{\ell(\nu)=\ell(\lambda)}
2^{z(\bar{\nu})+N'(\lambda,\nu)}\sigma_{\iota(\bar{\nu})}
+q^{-\delta(\lambda)-1}\sum_{\ell(\nu)=\ell(\lambda)+1}
2^{z(\bar{\nu})+N'(\lambda,\nu)}\sigma_{\iota(\bar{\nu})}.
\end{eqnarray}
On the other hand
$\sigma_{\iota(\lambda)}*\sigma_n=\sigma_{n,\iota(\lambda)}$
and the quantum Pieri rule gives 
\begin{eqnarray}\label{9} \nonumber
\iota(\sigma_{\lambda})*\iota(\sigma_k) = &2^{z(\lambda)+z(k)}
\Big(q^{-\delta(\lambda)-1}\sum_{\ell(\alpha)=\ell(\lambda)}
2^{N'(\iota(\lambda),\alpha)}\sigma_{n,\alpha}+
q^{-\delta(\lambda)-1}\sum_{\ell(\alpha)=\ell(\lambda)+1}
2^{N'(\iota(\lambda),\alpha)}\sigma_{n,\alpha} \\
 &+q^{-\delta(\lambda)}\sum_{\ell(\beta)=\ell(\lambda)+1}
2^{N'((n,\iota(\lambda)),\beta)}\sigma_{\bar{\beta}}
+q^{-\delta(\lambda)}\sum_{\ell(\beta)=\ell(\lambda)+2}
2^{N'((n,\iota(\lambda)),\beta)}\sigma_{\bar{\beta}}\Big).
\end{eqnarray}
As above we check that the first, second, third and fourth partial 
sums in (\ref{9}) coincide respectively with the second, fourth, third
and first partial sum in (\ref{8}). 

\smallskip\noindent {\it Second subcase}: $\iota(\lambda)_1=n$. Same story!
\qed

\medskip
After renormalizing, we get the following statement:

\begin{theo}\label{inv-OG}
The correspondence 
$$q\mapsto \frac{4}{q}, \hspace{15mm}
\sigma_{\lambda} \mapsto \Big(\frac{2}{q} \Big)^{\delta(\lambda)}
2^{\delta_{\iota(\lambda)_1,n}-\delta_{\lambda_1,n}}\sigma_{\iota(\lambda)} 
\quad for\;\lambda\subset\rho_n,$$
defines a ring involution of $QA^*(G_Q(n+1,2n+2))_{loc}$.
\end{theo}

\subsection{The Cayley plane}

For the Cayley plane $\OO\PP^2=E_6/P_1$, we have $d_{max}=2$, and 
the partition of the Hasse diagram is as follows. The two non trivial
intervals are isomorphic with the Hasse diagrams of an orthogonal
Grassmannian $G_Q(5,10)$, and a quadric $\QQ^8$.
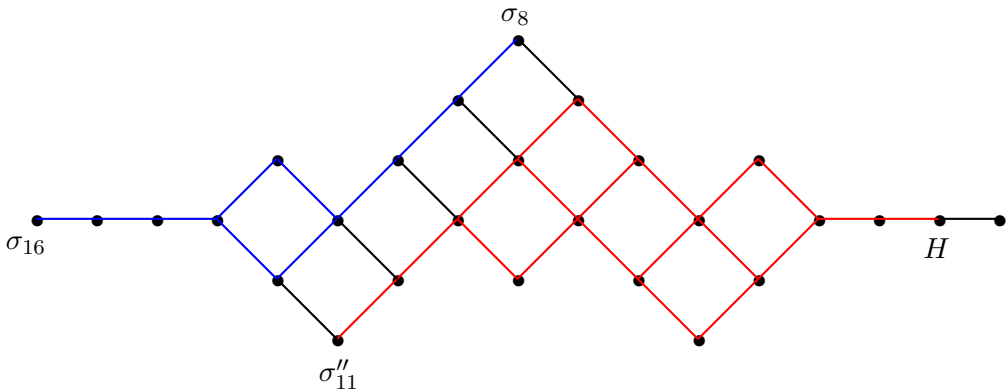
\begin{center}
\psset{unit=4mm}
\psset{xunit=4mm}
\psset{yunit=4mm}
\begin{pspicture*}(40,16)(-10,2)
\multiput(-6.2,7.7)(2,0){4}{$\bullet$}
\multiput(1.8,9.7)(4,0){5}{$\bullet$}
\multiput(1.8,5.7)(4,0){5}{$\bullet$}
\multiput(3.8,7.7)(4,0){5}{$\bullet$}
\multiput(3.8,3.7)(12,0){2}{$\bullet$}
\multiput(7.8,11.7)(4,0){2}{$\bullet$}
\multiput(21.8,7.7)(2,0){3}{$\bullet$}
\put(9.8,13.7){$\bullet$}
\put(-7,7){$\sigma_{16}$}
\put(9.45,14.6){$\s_8$}
\put(23.5,6.7){$H$}
\put(3.4,2.6){$\s_{11}''$}
\psline[linecolor=blue](-6,8)(0,8)
\psline[linecolor=blue](0,8)(2,10)
\psline[linecolor=blue](0,8)(2,6)
\psline[linecolor=blue](2,6)(10,14)
\psline[linecolor=red](4,4)(12,12)
\psline(10,14)(12,12)
\psline[linecolor=red](12,12)(18,6)
\psline[linecolor=red](10,6)(14,10)
\psline(2,6)(4,4)
\psline[linecolor=blue](2,10)(4,8)\psline(4,8)(6,6)
\psline(6,10)(8,8)
\psline(8,12)(10,10)
\psline[linecolor=red](10,10)(16,4)
\psline[linecolor=red](8,8)(10,6)
\psline[linecolor=red](16,4)(20,8)
\psline[linecolor=red](14,6)(18,10)
\psline[linecolor=red](18,10)(20,8)
\psline[linecolor=red](20,8)(24,8)
\psline(24,8)(26,8)
\end{pspicture*}
\end{center}

To each Schubert class $\sigma_w$ of $\OO\PP^2$ we associate a 
coefficient $\zeta(w)$ as follows, where $y=2x^2$ and $3y^2=1$.

\begin{center}
\psset{unit=4mm}
\psset{xunit=4mm}
\psset{yunit=4mm}
\begin{pspicture*}(40,15)(-10,2)
\put(-6.2,7.7){$1$}\put(-4.2,7.7){$\frac{1}{x}$}
\put(-2.2,7.7){$\frac{1}{y}$}\put(-.2,7.7){$2x$}
\put(1.8,9.7){$1$}\put(5.8,9.7){$y$}\put(9.8,9.7){$1$}
\put(13.8,9.7){$\frac{1}{y}$}\put(17.8,9.7){$1$}
\put(1.8,5.7){$1$}\put(5.8,5.7){$\frac{1}{y}$}\put(9.8,5.7){$1$}
\put(13.8,5.7){$y$}\put(17.8,5.7){$1$}
\put(3.8,7.7){$\frac{1}{2x}$}\put(7.8,7.7){$2x$}
\put(11.8,7.7){$\frac{1}{2x}$}\put(15.8,7.7){$2x$}
\put(17.8,7.7){$\frac{1}{2x}$}
\put(3.8,3.7){$\frac{1}{x}$}\put(15.8,3.7){$x$}
\put(7.8,11.7){$x$}\put(11.8,11.7){$\frac{1}{x}$}
\put(21.8,7.7){$\frac{1}{y}$}\put(23.8,7.7){$x$}
\put(25.8,7.7){$1$}
\put(9.8,13.7){$1$}
\end{pspicture*}
\end{center}

Observe that $\zeta(\iota(w))=\zeta(w^*)=\zeta(w)^{-1}$. 

\begin{prop} 
Theorem \ref{qinv} holds for the Cayley plane with 
$\kappa=12^{-\frac{1}{4}}$ and $\zeta(w)$ as above. 
\end{prop}

\proof A presentation for the quantum Chow ring $QA^*(\OO\PP^2)$ has 
been given in \cite{cmp}, Theorem 5.1, with generators $H$ and 
$\sigma'_4$. What we have to check is that the same relations are 
verified by their images 
$$\iota(H)=q^{-1}\sigma''_{11} \qquad {\rm and}\qquad 
\iota(\sigma'_4)=q^{-1}\sigma''_8.$$
This is a lengthy but direct computation.\qed
 
\subsection{The Freudenthal variety}

In this case $d_{max}=3$. The partition of the Hasse diagram is
symmetric, with two intervals that reduce to a point, and two 
that are isomorphic with the Hasse diagram of the Cayley plane.

\begin{center}
\psset{unit=3mm}
\psset{xunit=3mm}
\psset{yunit=3mm}
\begin{pspicture*}(50,20)(-10,-3)
\psline(-8,8)(-6,8)
\psline[linecolor=blue](-6,8)(0,8)
\psline[linecolor=blue](0,8)(4,12)
\psline[linecolor=blue](0,8)(2,6)
\psline[linecolor=blue](2,10)(10,2)\psline(10,2)(12,0)
\psline[linecolor=blue](2,6)(6,10)
\psline[linecolor=blue](6,6)(10,10)
\psline[linecolor=blue](4,12)(12,4)\psline(12,4)(14,2)
\psline[linecolor=blue](10,10)(14,6)\psline(14,6)(16,4)
\psline[linecolor=blue](16,12)(18,10)\psline(18,10)(20,8)
\psline[linecolor=blue](14,10)(16,8)\psline(16,8)(18,6)
\psline[linecolor=blue](8,4)(16,12)
\psline[linecolor=blue](10,2)(18,10)
\psline[linecolor=red](12,0)(20,8)
\psline[linecolor=blue](16,8)(18,8)
\psline[linecolor=red](18,6)(20,6)
\psline[linecolor=blue](18,10)(20,10)
\psline[linecolor=red](20,8)(22,8)
\psline[linecolor=blue](18,8)(26,16)
\psline(18,8)(20,6)\psline[linecolor=red](20,6)(22,4)
\psline[linecolor=red](20,6)(28,14)
\psline[linecolor=red](22,4)(30,12)
\psline(20,10)(22,8)\psline[linecolor=red](22,8)(24,6)
\psline(22,12)(24,10)\psline[linecolor=red](24,10)(28,6)
\psline(24,14)(26,12)\psline[linecolor=red](26,12)(34,4)
\psline(26,16)(28,14)\psline[linecolor=red](28,14)(36,6)
\psline[linecolor=red](28,6)(32,10)
\psline[linecolor=red](32,6)(36,10)
\psline[linecolor=red](34,4)(38,8)
\psline[linecolor=red](36,10)(38,8)
\psline[linecolor=red](38,8)(44,8)\psline(44,8)(46,8)
\multiput(-8.2,7.7)(2,0){5}{$\bullet$}
\multiput(1.7,9.7)(4,0){5}{$\bullet$}
\multiput(1.7,5.7)(4,0){5}{$\bullet$}
\multiput(3.7,7.7)(4,0){5}{$\bullet$}
\multiput(7.7,3.7)(4,0){3}{$\bullet$}
\multiput(3.7,11.7)(12,0){2}{$\bullet$}
\multiput(9.7,1.7)(4,0){2}{$\bullet$}
\put(11.7,-.3){$\bullet$}
\multiput(17.7,7.7)(4,0){5}{$\bullet$}
\multiput(19.7,9.7)(4,0){5}{$\bullet$}
\multiput(19.7,5.7)(4,0){5}{$\bullet$}
\multiput(21.7,11.7)(4,0){3}{$\bullet$}
\multiput(23.7,13.7)(4,0){2}{$\bullet$}
\multiput(25.7,15.7)(4,0){1}{$\bullet$}
\multiput(21.7,3.7)(12,0){2}{$\bullet$}
\multiput(37.7,7.7)(2,0){5}{$\bullet$}
\put(35.6,10.7){$\sigma'_5$}
\put(28.3,14.5){$\sigma_9$}
\put(28.3,5){$\sigma''_9$}
\put(11.9,-1.2){$\sigma_{17}$}
\put(43.5,8.5){$H$}
\end{pspicture*}
\end{center}

\medskip
In this case we needed a computer to check that Theorem \ref{SD} 
does hold. We know from \cite{cmp} that the quantum cohomology 
ring is generated over $\ZZ[q]$ by the classes $H,\sigma'_5$ and 
$\sigma_9$, and we know the relations explicitely. With the notations
of \cite{cmp}, our involution 
is given by 
$$q\mapsto 11943936q^{-1}, \quad H\mapsto q^{-1}\sigma_{17},
\quad \sigma'_5\mapsto 48q^{-1}\sigma'_{13},
\quad \sigma_9\mapsto 3456q^{-1}\sigma''_{9}.$$
We first checked that this map preserves the relations and is
involutive. Then 
we computed the image of each Schubert class and checked that 
it is given by the explicit form of Theorem \ref{SD}.

\section{Symmetries of Gromov-Witten invariants}

For convenience, denote by $p$ the symmetry of $W_X$ 
given by Poincar{\'e} duality. We deduce from Theorem \ref{qinv} 
the following identity for the quantum product:

\begin{theo}\label{sym0}
For any $u,v\in W_X$, we have
$$\sigma(\iota(u))*\sigma(p(v))=q^{\delta(u)+\delta(p(u))
-\delta(v)-\delta(p\iota(v))}\sigma(p\iota p(u))*
\sigma(\iota p\iota(v)).$$
\end{theo}

\proof Theorem \ref{qinv} is equivalent to the following identity
for Gromov-Witten invariants:
\begin{eqnarray}\label{GW1}
I_k(u,v,w)=\zeta(u)\zeta(v)\zeta(w)
I_{\delta(u)+\delta(v)-\delta(p(w))-k}(\iota(u),\iota(v),p\iota p(w)).
\end{eqnarray}
Using the fact that $\zeta(\iota(u))=\zeta(p(u))=\zeta(u)^{-1}$, the
same identity gives
$$I_{\ell}
(p\iota p(w),\iota(u),\iota(v))
=(\zeta(u)\zeta(v)\zeta(w))^{-1}
I_{\delta(p\iota p(w))+\delta(u)-\delta(p\iota(v))-\ell}
(\iota p\iota p(w),u,p\iota p\iota (v)).$$
Combining these two relations, we get
$$I_k(u,v,w)=I_{k+\delta(p(w))+\delta(p\iota
  p(w))-\delta(p\iota(v))-\delta(v)}
(u,p\iota p\iota (v),\iota p\iota p(w)),$$
which is equivalent to the identity for the quantum product
$$\sigma(p\iota p\iota (v))*\sigma(\iota p\iota p(w))=
q^{\delta(p(w))+\delta(p\iota
  p(w))-\delta(p\iota(v))-\delta(v)}\sigma(v)*\sigma(w).$$
Replacing $v$ by $\iota(u)$ and $w$ by $p(v)$ yields our claim. \qed

\medskip
Note that Theorem \ref{sym0} is non trivial only when $p$ and 
$\iota$ don't commute.

\smallskip
Now we observe that Theorem \ref{point}
can be used to generate more symmetries for the Gromov-Witten
invariants. 
The following statement is
a generalization of Proposition 4.10 in \cite{cmp}.

\begin{coro}\label{max}
For any $u\in W_X$, we have
\begin{eqnarray*}
 \; [Y_{d_{max}}] *\sigma(u) &= &q^{\delta(u)+\delta(p\iota(u))-d_{max}}
\sigma(p\iota p\iota(u)), \\
 \; [Y_{d_{max}}^*] *\sigma(u) &= &q^{d_{max}-\delta(p(u))}
\sigma(\iota p(u)).
\end{eqnarray*} 
\end{coro}

\proof To prove the first identity, 
multiply the identity of Theorem \ref{point} by $\sigma(pt)$
and use the fact that
$\sigma(pt)*\sigma(pt)=q^{d_{max}}[Y_{d_{max}}]$.
For the second identity, observe that Theorem \ref{point} implies 
that $\sigma(pt)$ is invertible in $QA^*(X)_{loc}$, and its inverse
verifies the formula
$$\sigma(pt)^{-1}*\sigma(u)=q^{-\delta(p(u))}\sigma(\iota p(u)).$$
Applying this to the fundamental class $\sigma(1)=1$ yields 
$\sigma(pt)^{-1}=q^{-d_{max}}[Y_{d_{max}}^*]$, which we just need to 
substitute in the previous identity.
\qed 

\begin{coro}\label{sym1}
For any $u,v\in W_X$, we have
$$q^{\delta(u)}\sigma(p(u))*\sigma(\iota(v))=
q^{\delta(v)}\sigma(\iota(u))*\sigma(p(v)).$$
\end{coro}

\proof Multiply the identity of Theorem \ref{point} for $\iota(u)$ by 
the Schubert class $\sigma(\iota(v))$, and use the associativity of the 
quantum product. \qed

\medskip Together with Theorem \ref{qinv}, we get a series of symmetry
relations for the Gromov-Witten invariants, which are generated by the 
following simple ones:

\begin{coro}\label{symGW}
For any $u,v,w\in W_X$, we have the relation
$$I_k(u,v,w)=\zeta(u)\zeta(v)\zeta(w)
I_{\delta(w)-k}(p(u),p(v),\iota(w)).$$
\end{coro}

\proof In terms of Gromov-Witten invariants, Corollary \ref{sym1}
writes
$$I_{k-\delta(u)}(p(u),\iota(v),w)=
I_{k-\delta(v)}(\iota(u),p(v),w).$$
Combining with the identity (\ref{GW1}), we deduce that 
\begin{eqnarray*}
I_k(u,v,w) &= &I_{k+\delta(p(u))-\delta(v)}(\iota p(u),p \iota(v),w) \\
  &= &I_{k+\delta(p(u))-\delta(v)}(\iota p(u),w, p \iota(v)) \\
 & =&\zeta(\iota p(u))\zeta(w)\zeta(p \iota(v))
I_{\delta(p(u))+\delta(w)-\delta(v))-(k+\delta(p(u))-\delta(v))}
(p(u),\iota(w),p(v))\\
 & =&\zeta(u)\zeta(w)\zeta(v)I_{\delta(w)-k}(p(u),p(v),\iota(w)).
\end{eqnarray*}
This is what we wanted to prove. \qed

\begin{coro}\label{qmax}
For $u,v\in W_X$, the maximal power of $q$ that appears in the 
quantum product of the Schubert classes $\sigma(u)*\sigma(v)$ is
$$d_{max}(u,v)=\delta(u)-\delta(\iota(u),p(v))
=\delta(v)-\delta(\iota(v),p(u)).$$
\end{coro}

Corollary \ref{symGW} suggests to study the group $\Gamma $ of permutations of 
$W_X^3$ generated by $(p,p,\iota)$, $(p,\iota,p)$, $(\iota,p,p)$. 
Clearly the size of this group is governed by the order $\eta$ of the 
permutation $p\iota$ of $W_X$. 

\begin{prop}\label{order} The order of $\Gamma$ is $2\eta^2$, 
and $\eta$ is given by the following table:
$$\begin{array}{ccc}
X & \qquad & \eta \\
 & & \\
\G(p,n) & & n/gcd(p,n-p) \\
\G_{\omega}(n,2n) & & 2\\
\G_Q(n,2n) & & 4/gcd(2,n)\\
\OO\PP^2 & & 3\\
E_7/P_7 & & 2 
\end{array}$$ 
\end{prop}
 
Note that $\eta$ always divides the order of the symmetry group
of the affine Dynkin diagram of $G$. 
\medskip

\proof Let $\Gamma_0$ denote the group of permutations of $W_X$
generated by the two involutions $p$ and $\iota$. The order of 
$\Gamma_0$ is $2\eta$. Moreover, the projection on the first factor 
yields a morphism $\Gamma\ra\Gamma_0$ which is obviously surjective.
Its kernel consists in the permutations of type $(1,(pi)^k,(ip)^k)$, 
with $k\in\ZZ$, so its order is $\eta$. Thus the order of $\Gamma$ is 
$2\eta^2$. 

For the explicit values of $\eta$, first consider the case of
$X=\G(p,n)$. A partition $\lambda\in W_X$ can be identified with 
a $01$-sequence $\omega$ with $p$ ones and $n-p$ zeroes encoding vertical
and horizontal steps along the boundary of $\lambda$, starting 
from the SW corner. Then the size of the biggest square contained in 
$\lambda$ is the number of zeroes among the first $p$ terms of the sequence. 
Moreover, reading the sequence backwards we get
the sequence $\omega^*$ of the Poincar{\'e} dual partition $p(\lambda)$. 
So to get the sequence $\omega'$ of $\iota(\lambda)$, we write $\omega=
\omega_0\omega_1$ where $\omega_0$ has length $p$ and $\omega_1$
has length $n-p$, and let $\omega'=\omega_0^*\omega_1^*$. To 
deduce the  sequence $\omega''$ of $p\iota(\lambda)$ we simply reverse 
$\omega'$, so $\omega''=\omega_1\omega_0$. The claim easily follows. 

Now suppose $X=\G_Q(n,2n)$. We identify a strict partition 
$\lambda\in W_X$ with a $01$-sequence $\omega$ of length $n$
as follows. First we consider it as a usual partition in a square of
size $(n+1)\times (n+1)$ and we let $\omega'$ be the associated
$01$-sequence, of length $2n+2$. It begins with a $1$ and ends with a
$0$. We suppress the initial $1...10$ sequence. Moreover, since 
$\lambda$ is strict every $1$ is followed by a $0$, which we 
suppress. The resulting sequence has length $n$ and is our 
$\omega=\omega_1\cdots\omega_n$. Note that the length of $\lambda$ is
the number of $1$'s.
 
We check that $p$ and $\iota$ are easily expressed in terms of 
$01$-sequences:
$$p(\lambda)\mapsto \bar{\omega}_1\cdots \bar{\omega}_n, \hspace{2cm}
\iota(\lambda)\mapsto \omega_{n-1}\cdots\omega_1\omega_0,$$
 where $\omega_0=\omega_1+\cdots +\omega_n\; (mod\; 2)$ and 
$\bar{0}=1$, $\bar{1}=0$. So $p\iota(\lambda)\mapsto
\bar{\omega}_{n-1}\cdots \bar{\omega}_1\bar{\omega}_0$, and 
$p\iota p\iota(\lambda)\mapsto \omega_1\cdots\omega_{n-1}\omega'_n$, 
where $\omega'_n=\omega_n$ if $n$ is even and
$\omega'_n=\bar{\omega}_n$ if $n$ is odd.

The case of $\G_{\omega}(n,2n)$, and also that of $E_7/P_7$ are trivial, 
since $p$ commutes with $\iota$.  Finally the case of $\OO\PP^2$ follows
from an explicit computation. \qed

\medskip For a Grassmannian $\G(p,n)$, we get a $\ZZ_n$-symmetry only
when $(p,n)=1$, while this symmetry always exists by \cite{Po}. 
Does the same phenomenon happen for $\G_Q(2m,4m)$ ? That is, does our
twofold symmetry extend to a fourfold symmetry ? 

\medskip The cases for which $\eta=2$ are the most symmetric: the
involution $\iota$ commutes with Poincar{\'e} duality. In particular
we  get the relation
$$\delta(u)+\delta(p(u))=d_{max}.$$
The Gromov-Witten invariants are identified by groups of eight
according to the following identities:
\begin{eqnarray*}
I_k(u,v,w) &= &\zeta(u)\zeta(v)\zeta(w)
I_{\delta(w)-k}(p(u),p(v),\iota(w)) \\
 &= & I_{k+d_{\max}-\delta(v)-\delta(w)}(u,\iota p(v),\iota p(w)) \\
&= &\zeta(u)\zeta(v)\zeta(w)
I_{\delta(u)+\delta(v)+\delta(w)-d_{max}-k}(\iota(u),\iota(v),\iota(w)). 
\end{eqnarray*}
In particular all Gromov-Witten invariants can be directly computed 
from those of degree $k\le d_{max}/4$. 
 
\medskip
For the other cases we get even more identities, so that lots of 
Gromov-Witten invariants can be identified with classical intersection
numbers on the same variety. 

\medskip 
We close this paper with a ``dual quantum Chevalley formula'', 
that we obtain by applying our strange duality Theorem to 
the quantum Chevalley formula of \cite{cmp}. Remember that
$\alpha_0$ denotes the highest root, and that
$\sigma(s_{\alpha_0})=[T_1]$. The Bruhat interval $[1,T_1]$ has 
a Poincar{\'e} involution $p_{T_1}$. Finally, recall that we denoted 
by $\beta$  the simple root that defines $P$.   

\begin{prop}
For any $u\in W_X$, we have
$$\sigma(u)*\sigma(s_{\alpha_0})=\delta_{p(u)\le\alpha_0}
\sigma(p_{T_1}p(u))+q\sum_{\substack{s_{\alpha}u\rightarrow u,\\ \alpha\ne\beta}}
n_{\alpha}(\alpha_0)\sigma(s_{\alpha}u)
+q\delta_{\iota(u)\le\alpha_0}y(p(s_{\alpha_0}))\sigma(s_{\beta}u).$$
\end{prop}

The classical intersection product $\delta_{p(u)\le\alpha_0}\sigma(p_{T_1}p(u))$
is in agreement with Proposition \ref{prolisse}. It would be
interesting to extend this formula to more general products.

\bigskip\noindent
Pierre-Emmanuel {\sc Chaput}, Laboratoire de Math{\'e}matiques Jean Leray, 
UMR 6629 du CNRS, UFR Sciences et Techniques,  2 rue de la Houssini{\`e}re, BP
92208, 44322 Nantes cedex 03, France. 

\noindent {\it email}: pierre-emmanuel.chaput@math.univ-nantes.fr

\medskip\noindent
Laurent {\sc Manivel}, 
Institut Fourier, UMR 5582 du CNRS,  Universit{\'e} de Grenoble I, 
BP 74, 38402 Saint-Martin d'H{\'e}res, France.
 
\noindent {\it email}: Laurent.Manivel@ujf-grenoble.fr

\medskip\noindent
Nicolas {\sc Perrin},  Institut de Math{\'e}matiques,  
Universit{\'e} Pierre et Marie Curie, Case 247, 4 place Jussieu,  
75252 PARIS Cedex 05, France.

\noindent {\it email}: nperrin@math.jussieu.fr
\end{document}